\documentclass[a4paper,11pt]{amsart}

\usepackage{graphicx}
\usepackage{mathptmx}
\usepackage{amsmath}
\usepackage{amssymb}
\usepackage{enumitem}
\usepackage{xcolor}

\newmuskip\pFqmuskip

\newcommand*\pFq[6][8]{%
  \begingroup 
  \pFqmuskip=#1mu\relax
  \mathcode`=\string"8000
  \begingroup\lccode`\~=`\,
  \lowercase{\endgroup\let~}\pFqcomma
  F^{#2}_{#3}{\left(\genfrac..{0pt}{}{#4}{#5}\bigg|#6\right)}%
  \endgroup
}
\newcommand{\pFqcomma}{\mskip\pFqmuskip}

\newtheorem{theorem}{Theorem}[section]

\begin{document}

\title[]{Probabilistic degenerate $r$-Stirling numbers of the  second and probabilistic degenerate $r$-Bell polynomials}

\author{Taekyun  Kim}
\address{Department of Mathematics, Kwangwoon University, Seoul 139-701, Republic of Korea}
\email{tkkim@kw.ac.kr}
\author{Dae San  Kim}
\address{Department of Mathematics, Sogang University, Seoul 121-742, Republic of Korea}
\email{dskim@sogang.ac.kr}


\begin{abstract}
Assume that $Y$ is a random variable whose moment generating function exists in a neighborhood of the origin. We study the probabilistic degenerate $r$-Stirling numbers of the second kind associated with $Y$ and the probabilistic degenerate $r$-Bell polynomials associated with $Y$. They are respectively probabilistic extensions of the degenerate $r$-Stirling numbers of the second and the degenerate $r$-Bell polynomials which are studied earlier. The aim of this paper is to obtain their properties, explicit expressions and related identities by means of generating functions.
\end{abstract}

\maketitle

\markboth{\centerline{\scriptsize Probabilistic degenerate $r$-Stirling numbers of the  second and probabilistic degenerate $r$-Bell polynomials}}
{\centerline{\scriptsize T. Kim and D. S. Kim}}

\section{Introduction}
Let $Y$ be a random variable whose moment generating function exists in a neighborhood of the origin. In this paper, we study the probabilistic extensions of the degenerate $r$-Stirling numbers of the second kind and the degenerate $r$-Bell polynomials. Namely, they are respectively the probabilistic degenerate $r$-Stirling numbers of the second kind associated with $Y$ and the probabilistic degenerate $r$-Bell polynomials associated with $Y$. The aim of this paper is to obtain their properties, explicit expressions and related identities by means of generating functions. We remark here that degenerate versions (see  [12,16,18,23-27,30-32] and the references therein), $\lambda$-analogues (see [15,19,28] and the references therein) and probabilistic extensions (see [3,10,20-22,29,35] and the references therein) of many special numbers and polynomials have been intensively studied in recent years and a lot of interesting results have been obtained.\par
The outline of this paper is as follows. In Section 1, we recall the degenerate exponentials, the Stirling numbers of both kinds, the degenerate Stirling numbers of the second kind and the degenerate $r$-Stirling numbers of the second kind. We recall the probabilistic degenerate Stirling numbers of the second kind associated with $Y$, $S_{2,\lambda}^{Y}(n,k)$, and the probabilistic degenerate Bell polynomials associated with $Y$, $\mathrm{Bel}_{n,\lambda}^{Y}(x)$. In addition, their explicit expressions are given. Section 2 contains the main results of this paper. Let $Y$ be as before, and let $S_{0}=0$,\quad $S_{k}=Y_{1}+Y_{2}+\cdots+Y_{k},\ (k\ge 1)$. Let $r$ be a nonnegative integer. We define the probabilistic degenerate $r$-Stirling numbers of the second kind associated with $Y$, $S_{2,\lambda}^{(r,Y)}(n+r,k+r)$. We obtain three explicit expressions of $S_{2,\lambda}^{(r,Y)}(n+r,k+r)$ in Theorems 2.1-2.3. Indeed, we represent $S_{2,\lambda}^{(r,Y)}(n+r,k+r)$ as a finite sum involving $E[(S_{j+r})_{n,\lambda}]$ in Theorem 2.1, a finite sum involving $S_{1}(n-l,m)$, $S_{2,\lambda}^{Y}(l,k)$, and $E[S_{r}^{m}]$ in Theorem 2.2 and a finite sum involving $S_{2,\lambda}^{Y}(n,m+k)$ in Theorem 2.3. In Theorem 2.4, we derive a finite sum identity relating $S_{2,\lambda}^{(r,Y)}(n+r,k+r)$ to
$S_{2,\lambda}^{Y}(n,k)$. We define the probabilistic degenerate $r$-Bell polynomials associated with $Y$, $\mathrm{Bel}_{n,\lambda}^{(r,Y)}(x)$. We deduce three different expressions of $\mathrm{Bel}_{n,\lambda}^{(r,Y)}(x)$ in Theorems 2.5-2.7. Indeed, we show $\mathrm{Bel}_{n,\lambda}^{(r,Y)}(x)$ is equal to $\sum_{k=0}^{n}S_{2,\lambda}^{(r,Y)}(n+r,k+r)x^{k}$ in Theorem 2.5, represent it as a finite sum involving $E[(S_{r})_{m,\lambda}]$ and $\mathrm{Bel}_{n-m,\lambda}^{Y}(x)$ in Theorem 2.6 and derive for it a Dobinski-like formula in Theorem 2.7. We get a recurrence relation for $S_{2,\lambda}^{(r,Y)}(n+r,k+r)$ in Theorem 2.8 and a finite sum identity relating it to $S_{1}(j,i)$ and $S_{2,\lambda}^{Y}(n,i)$ in Theorem 2.9. As general references on polynomials, functions, combinatorics and probability, one may refer to [2,5,11,13,33,34]. In addition, the reader refers to [1,4,6-9,14,17] for some related special numbers and polynomials. For the rest of this section, we recall the facts that are needed throughout this paper.

\vspace{0.1in}

Let $[n]=\{1,2,3,\dots,n\}$. The Stirling number of the second kind $S_{2}(n,k)$ counts the number of partitions of the set $[n]$ into non-empty disjoint subsets. Let $r$ be a nonnegative integer. The $r$-Stirling number of the second kind $S_{2}^{(r)}(n,k)$ counts the number of partitions of the set $[n]$ into $k$ non-empty disjoint subsets in such a way that the numbers $1,2,\dots,r$ are in distinct subsets (see [1-35]). \par
For any nonzero $\lambda\in\mathbb{R}$, the degenerate exponentials are given by
\begin{equation}
e_{\lambda}^{x}(t)=(1+ \lambda t)^{\frac{x}{\lambda}}=\sum_{k=0}^{\infty}(x)_{k,\lambda}\frac{t^{k}}{k!},\quad e_{\lambda}(t)=e_{\lambda}^{1}(t),\quad (\mathrm{see}\ [14-32]), \label{1}
\end{equation}
where $(x)_{0,\lambda}=1,\ (x)_{n,\lambda}=x(x-\lambda)\cdots(x-(n-1)\lambda),\ (n\ge 1)$.\par
For $n\ge 0$, the Stirling numbers of the first kind are defined by
\begin{equation}
(x)_{n}=\sum_{k=0}^{n}S_{1}(n,k)x^{k},\quad (\mathrm{see}\ [1-30]), \label{2}	
\end{equation}
 where $(x)_{0}=1,\ (x)_{n}=x(x-1)(x-2)\cdots(x-n+1),\ (n\ge 1)$. \\
Alternatively, they are given by
\begin{equation}
\frac{1}{k!}\big(\log (1+t) \big)^{k}=\sum_{n=k}^{\infty}S_{1}(n,k)\frac{t^{n}}{n!}. \label{2-1}
\end{equation}
As the inversion of \eqref{2}, the Stirling numbers of the second kind are given by
 \begin{equation}
 x^{n}=\sum_{k=0}^{n}S_{2}(n,k)(x)_{k},\quad (n\ge 0),\quad (\mathrm{see}\ [14,30-35]). \label{3}
 \end{equation}
In [16], the degenerate Stirling numbers of the second kind are defined by
\begin{equation}
(x)_{n,\lambda}=\sum_{k=0}^{n}S_{2,\lambda}(n,k)(x)_{k},\quad (n\ge 0). \label{4}	
\end{equation}
Note that $\lim_{\lambda\rightarrow 0}S_{2,\lambda}(n,k)=S_{2}(n,k),\ (n\ge k\ge 0)$. \par
Let $r$ be a nonnegative integer. Then the degenerate $r$-Stirling numbers of the second kind are defined by
\begin{equation}
(x+r)_{n,\lambda}=\sum_{k=0}^{n}S_{2,\lambda}^{(r)}(n+r,k+r)(x)_{k},\quad (\mathrm{see}\ [30]).\label{5}
\end{equation}
Note that $\lim_{\lambda\rightarrow 0}S_{2,\lambda}^{(r)}(n+r,k+r)= S_{2}^{(r)}(n+r,k+r)$ are the $r$-Stirling numbers of the second kind given by
\begin{equation*}
(x+r)^{n}=\sum_{k=0}^{n}S_{2}^{(r)}(n+r,k+r)(x)_{k},\quad (n\ge 0).	
\end{equation*}
From \eqref{5}, we note that
\begin{equation}
\frac{1}{k!}\Big(e_{\lambda}(t)-1\Big)^{k}e_{\lambda}^{r}(t)=\sum_{n=k}^{\infty}S_{2,\lambda}^{(r)}(n+r,k+r)\frac{t^{n}}{n!},\quad (k\ge 0),\quad (\mathrm{see}\ [30]).\label{6}	
\end{equation}
Assume that $Y$ is a random variable such that the moment generating function of $Y$,
\begin{equation}
E\big[e^{tY}\big]=\sum_{n=0}^{\infty}E[Y^{n}]\frac{t^{n}}{n!},\quad (|t|<r) \label{6-1}
\end{equation}
exists for some $r>0$. Let $(Y_{j})_{j\ge 1}$ be a sequence of mutually independent copies of the random variable $Y$, and let
\begin{displaymath}
S_{k}=Y_{1}+Y_{2}+\cdots+Y_{k},\quad (k\ge 1),\ \mathrm{and}\ S_{0}=0.
\end{displaymath}
The probabilistic degenerate Stirling numbers of the second kind associated with $Y$ are defined by
\begin{equation}
\frac{1}{k!}\Big(E\big[e_{\lambda}^{Y}(t)\big]-1\Big)^{k}=\sum_{n=k}^{\infty}S_{2,\lambda}^{Y}(n,k)\frac{t^{n}}{n!},\quad (n\ge 0),\quad (\mathrm{see}\ [21,29]). \label{7}
\end{equation}
When $Y=1$, we have $S_{2,\lambda}^{Y}(n,k)=S_{2,\lambda}(n,k)$. \par
In addition, the probabilistic degenerate Bell polynomials associated with $Y$ are given by
\begin{equation}
e^{x(E[e_{\lambda}^{Y}(t)]-1)}=\sum_{n=0}^{\infty}\mathrm{Bel}_{n,\lambda}^{Y}(x)\frac{t^{n}}{n!},\quad (\mathrm{see}\ [21]). \label{8}	
\end{equation}
When $Y=1$, $\mathrm{Bel}_{n,\lambda}^{Y}(x)=\mathrm{Bel}_{n,\lambda}(x),\ (n\ge 0)$, where $\mathrm{Bel}_{n,\lambda}(x)$ are the degenerate Bell polynomials defined by
\begin{equation*}
\mathrm{Bel}_{n,\lambda}(x)=\sum_{k=0}^{n}S_{2,\lambda}(n,k)x^{k},\quad (n\ge 0),\quad (\mathrm{see}\ [20,23,30]).
\end{equation*}
In particular, for $x=1,\ \mathrm{Bel}_{n,\lambda}^{Y}=\mathrm{Bel}_{n,\lambda}^{Y}(1)$ are called the probabilistic degenerate Bell numbers. \par
From \eqref{7} and \eqref{8}, we note that
\begin{equation}
S_{2,\lambda}^{Y}(n,k)=\frac{1}{k!}\sum_{j=0}^{k}\binom{k}{j}(-1)^{k-j}E\big[(S_{j})_{n,\lambda}\big],\quad (n\ge k\ge 0), \label{9}	
\end{equation}
and
\begin{equation}
\mathrm{Bel}_{n,\lambda}^{Y}(x)=\sum_{k=0}^{n}S_{2,\lambda}^{Y}(n,k)x^{k},\quad (n\ge 0),\quad (\mathrm{see}\ [21]).
\end{equation}

\section{Probabilistic degenerate $r$-Stirling numbers of the second kind and probabilistic degenerate $r$-Bell polynomials}
Let $(Y_{j})_{j\ge 1}$ be a sequence of mutually independent copies of the random variable $Y$, and let
\begin{displaymath}
S_{0}=0,\quad\mathrm{and}\quad S_{k}=Y_{1}+Y_{2}+\cdots+Y_{k},\ (k\ge 1).
\end{displaymath}
Let $r$ be a nonnegative integer. Now, we consider the {\it{probabilistic degenerate $r$-Stirling numbers of the second kind associated with $Y$}} given by
\begin{equation}
\frac{1}{k!}\Big(E\big[e_{\lambda}^{Y}(t)\big]-1\Big)^{k} E\big[e_{\lambda}^{S_{r}}(t)\big]=\sum_{n=k}^{\infty}S_{2,\lambda}^{(r,Y)}(n+r,k+r)\frac{t^{n}}{n!}.\label{11}
\end{equation}
When $Y=1$, $S_{2,\lambda}^{(r,Y)}(n+r,k+r)=S_{2,\lambda}^{(r)}(n+r,k+r),\ (n\ge k\ge 0)$, (see \eqref{5}).\par
From \eqref{11} and \eqref{1}, we note that
\begin{align}
\sum_{n=k}^{\infty}S_{2,\lambda}^{(r,Y)}(n+r,k+r)\frac{t^{n}}{n!}&=\frac{1}{k!}\Big(E\big[e_{\lambda}^{Y}(t)\big]-1\Big)^{k} E\big[e_{\lambda}^{S_{r}}(t)\big]\label{12} \\
&=\frac{1}{k!}\sum_{j=0}^{k}\binom{k}{j}(-1)^{k-j}\Big(E\big[e_{\lambda}^{Y}\big]\Big)^{j}\Big(E\big[e_{\lambda}^{Y}(t)\big]\Big)^{r}\nonumber \\
&=\frac{1}{k!}\sum_{j=0}^{k}\binom{k}{j}(-1)^{k-j}E\big[e_{\lambda}^{S_{j+r}}(t)\big]\nonumber \\
&=\sum_{n=0}^{\infty}\frac{1}{k!}\sum_{j=0}^{k}\binom{k}{j}(-1)^{k-j}E\big[(S_{j+r})_{n,\lambda}\big]\frac{t^{n}}{n!}.\nonumber	
\end{align}
By comparing the coefficients on both sides of \eqref{12}, we obtain the following theorem.
\begin{theorem}
For $n\ge k\ge 0$, we have
\begin{displaymath}
S_{2,\lambda}^{(r,Y)}(n+r,k+r)=\frac{1}{k!}\sum_{j=0}^{k}\binom{k}{j}(-1)^{k-j}E\big[(S_{j+r})_{n,\lambda}\big].
\end{displaymath}
\end{theorem}
By \eqref{1}, \eqref{2-1} and \eqref{7}, we get
\begin{align}
\frac{1}{k!}\Big(E\big[e_{\lambda}^{Y}(t)\big]-1\Big)^{k}&E\big[e_{\lambda}^{S_{r}}(t)\big]=\sum_{l=k}^{\infty}S_{2,\lambda}^{Y}(l,k)\frac{t^{l}}{l!}\sum_{m=0}^{\infty}\frac{E[S_{r}^{m}]}{\lambda^{m}}\frac{1}{m!}\big(\log(1+\lambda t)\big)^{m}\label{13} \\
&=\sum_{l=k}^{\infty}S_{2,\lambda}^{Y}(l,k)\frac{t^{l}}{l!}\sum_{m=0}^{\infty}E[S_{r}^{m}]\lambda^{-m}\sum_{i=m}^{\infty}S_{1}(i,m)\lambda^{i}\frac{t^{i}}{i!}\nonumber \\
&=\sum_{l=k}^{\infty}S_{2,\lambda}^{Y}(l,k)\frac{t^{l}}{l!}\sum_{i=0}^{\infty}\sum_{m=0}^{i}E[S_{r}^{m}]\lambda^{i-m}S_{1}(i,m)\frac{t^{i}}{i!}\nonumber \\
&=\sum_{n=k}^{\infty}\bigg(\sum_{l=k}^{n}\sum_{m=0}^{n-l}\binom{n}{l}\lambda^{n-m-l}E\big[S_{r}^{m}\big]S_{1}(n-l,m)S_{2,\lambda}^{Y}(l,k)\bigg)\frac{t^{n}}{n!}.\nonumber
\end{align}
Therefore, by \eqref{11} and \eqref{13}, we obtain the following theorem.
\begin{theorem}
For $n\ge k\ge 0$, we have
\begin{displaymath}
S_{2,\lambda}^{(r,Y)}(n+r,k+r)=\sum_{l=k}^{n}\sum_{m=0}^{n-l}\binom{n}{l}\lambda^{n-m-l}S_{1}(n-l,m)S_{2,\lambda}^{Y}(l,k)E\big[S_{r}^{m}\big].
\end{displaymath}
\end{theorem}
From \eqref{11} and \eqref{7}, we have
\begin{align}
\sum_{n=k}^{\infty}S_{2,\lambda}^{(r,Y)}(n+r,k+r)\frac{t^{n}}{n!}&=\frac{1}{k!}\Big(E\big[e_{\lambda}^{Y}(t)\big]-1\Big)^{k}E\big[e_{\lambda}^{S_{r}}(t)\big] \label{14} \\
&=\frac{1}{k!}\Big(E\big[e_{\lambda}^{Y}(t)\big]-1\Big)^{k}\Big(E\big[e_{\lambda}^{Y}(t)\big]-1+1\Big)^{r}\nonumber \\
&=\frac{1}{k!}\sum_{m=0}^{\infty}\binom{r}{m}\Big(E\big[e_{\lambda}^{Y}(t)\big]-1\Big)^{m+k} \nonumber \\
&=\sum_{m=0}^{\infty}\binom{r}{m}m!\binom{m+k}{m}\sum_{n=m+k}^{\infty}S_{2,\lambda}^{Y}(n,m+k)\frac{t^{n}}{n!}\nonumber\\
&=\sum_{n=k}^{\infty}\bigg(\sum_{m=0}^{n-k}\binom{r}{m}m!\binom{m+k}{m}S_{2,\lambda}^{Y}(n,m+k)\bigg)\frac{t^{n}}{n!}.\nonumber
\end{align}
Therefore, by comparing the coefficients on both sides of \eqref{4}, we obtain the following theorem.
\begin{theorem}
For $n,k\ge 0$ with $n\ge k$, we have
\begin{displaymath}
S_{2,\lambda}^{(r,Y)}(n+r,k+r)=\sum_{m=0}^{n-k}\binom{m+k}{m}\binom{r}{m}m!S_{2,\lambda}^{Y}(n,m+k).
\end{displaymath}
\end{theorem}
By binomial expansion and \eqref{11}, we get
\begin{align}
\Big(	E\big[e_{\lambda}^{Y}(t)\big]\Big)^{x+r}&=\Big(E\big[e_{\lambda}^{Y}(t)\big]\Big)^{r}\Big(E\big[e_{\lambda}^{Y}(t)\big]-1+1\Big)^{x}\label{15} \\
&=\sum_{k=0}^{\infty}(x)_{k}\frac{1}{k!}\Big(E\big[e_{\lambda}^{Y}(t)\big]-1\Big)^{k}\Big(E\big[e_{\lambda}^{Y}(t)\big]\Big)^{r}\nonumber \\
&=\sum_{k=0}^{\infty}(x)_{k}\sum_{n=k}^{\infty}S_{2,\lambda}^{(r,Y)}(n+r,k+r)\frac{t^{n}}{n!}\nonumber \\
&=\sum_{n=0}^{\infty}\bigg(\sum_{k=0}^{n}(x)_{k}S_{2,\lambda}^{(r,Y)}(n+r,k+r)\bigg)\frac{t^{n}}{n!},\nonumber
\end{align}
and, by \eqref{7}, we have
\begin{align}
\Big(E\big[e_{\lambda}^{Y}(t)\big]\Big)^{x+r}&=\sum_{k=0}^{\infty}(x+r)_{k}\frac{1}{k!}\Big(E\big[e_{\lambda}^{Y}(t)\big]-1\Big)^{k}\label{16}\\
&=\sum_{k=0}^{\infty}(x+r)_{k}\sum_{n=k}^{\infty}S_{2,\lambda}^{Y}(n,k)\frac{t^{n}}{n!}\nonumber \\
&=\sum_{n=0}^{\infty}\bigg(\sum_{k=0}^{n}(x+r)_{k}S_{2,\lambda}^{Y}(n,k)\bigg)\frac{t^{n}}{n!}.\nonumber
\end{align}
Therefore, by \eqref{15} and \eqref{16}, we obtain the following theorem.
\begin{theorem}
For $n \ge 0$, we have
\begin{displaymath}
\sum_{k=0}^{n}S_{2,\lambda}^{(r,Y)}(n+r,k+r)(x)_{k}=\sum_{k=0}^{n}S_{2,\lambda}^{Y}(n,k)(x+r)_{k}.
\end{displaymath}
\end{theorem}
Now, we consider the {\it{probabilistic degenerate $r$-Bell polynomials associated with $Y$}} given by
\begin{align}
E\big[e_{\lambda}^{S_{r}}(t)\big]e^{x(E[e_{\lambda}^{Y}(t)]-1)}=\sum_{n=0}^{\infty}\mathrm{Bel}_{n,\lambda}^{(r,Y)}(x)\frac{t^{n}}{n!}. \label{17}
\end{align}
In particular, for $x=1$, $\mathrm{Bel}_{n,\lambda}^{(r,Y)}=\mathrm{Bel}_{n,\lambda}^{(r,Y)}(1)$ are called the {\it{probabilistic degenerate $r$-Bell numbers associated with $Y$}}. When $Y=1$, $\mathrm{Bel}_{n,\lambda}^{(r,Y)}(x)=\mathrm{Bel}_{n,\lambda}^{(r)}(x)$ is the degenerate $r$-Bell polynomials given by (see [30])
\begin{equation}
e_{\lambda}^{r}(t)e^{x(e_{\lambda}(t)-1)}=\sum_{n=0}^{\infty}\mathrm{Bel}_{n,\lambda}^{(r)}(x)\frac{t^{n}}{n!}. \label{17-1}
\end{equation}
By \eqref{11}, we note that
\begin{align}
E\big[e_{\lambda}^{S_{r}}(t)\big]e^{x(E[e_{\lambda}^{Y}(t)]-1)}&=\sum_{k=0}^{\infty}x^{k}\frac{1}{k!}\Big(E\big[e_{\lambda}^{Y}(t)\big]-1\Big)^{k}E\big[e_{\lambda}^{S_{r}}(t)\big] \label{18} \\
&=\sum_{k=0}^{\infty}x^{k}\sum_{n=k}^{\infty}S_{2,\lambda}^{(r,Y)}(n+r,k+r)\frac{t^{n}}{n!}\nonumber\\
&=\sum_{n=0}^{\infty}\bigg(\sum_{k=0}^{n}x^{k}S_{2,\lambda}^{(r,Y)}(n+r,k+r)\bigg)\frac{t^{n}}{n!}. \nonumber
\end{align}
Therefore, by \eqref{17} and \eqref{18}, we obtain the following theorem.
\begin{theorem}
For $n\ge 0$, we have
\begin{displaymath}
\mathrm{Bel}_{n,\lambda}^{(r,Y)}(x)=\sum_{k=0}^{n}S_{2,\lambda}^{(r,Y)}(n+r,k+r)x^{k}.
\end{displaymath}
\end{theorem}
By \eqref{1}, \eqref{8} and \eqref{17}, we get
\begin{align}
\sum_{n=0}^{\infty}\mathrm{Bel}_{n,\lambda}^{(r,Y)}(x)\frac{t^{n}}{n!}&=E\big[e_{\lambda}^{S_{r}}(t)\big]e^{x(E[e_{\lambda}^{Y}(t)]-1)}\label{20} \\
&=\sum_{m=0}^{\infty}E\big[(S_{r})_{m,\lambda}\big]\frac{t^{m}}{m!}\sum_{l=0}^{\infty}\mathrm{Bel}_{l,\lambda}^{Y}(x)\frac{t^{l}}{l!}\nonumber \\
&=\sum_{n=0}^{\infty}\bigg(\sum_{m=0}^{n}\binom{n}{m}E\big[(S_{r})_{m,\lambda}\big]\mathrm{Bel}_{n-m,\lambda}^{Y}(x)\bigg)\frac{t^{n}}{n!}.\nonumber
\end{align}
Therefore, by \eqref{20}, we obtain the following theorem.
\begin{theorem}
For $n\ge 0$, we have
\begin{displaymath}
\mathrm{Bel}_{n,\lambda}^{(r,Y)}(x)=\sum_{m=0}^{n}\binom{n}{m}E\big[(S_{r})_{m,\lambda}\big]\mathrm{Bel}_{n-m,\lambda}^{Y}(x).
\end{displaymath}
\end{theorem}
Now, we observe from \eqref{1} and \eqref{17} that
\begin{align}
\sum_{n=0}^{\infty}\mathrm{Bel}_{n,\lambda}^{(r,Y)}(x)\frac{t^{n}}{n!}&=E\big[e_{\lambda}^{S_{r}}(t)\big]e^{x(E[e_{\lambda}^{Y}(t)]-1)}\label{21} \\
&=e^{-x}\sum_{k=0}^{\infty}\frac{x^{k}}{k!}\Big(E\big[e_{\lambda}^{Y}(t)\big]\Big)^{k+r}\nonumber \\
&=e^{-x}\sum_{k=0}^{\infty}\frac{x^{k}}{k!}\sum_{n=0}^{\infty}E\big[(S_{k+r})_{n,\lambda}\big]\frac{t^{n}}{n!}\nonumber \\
&=\sum_{n=0}^{\infty}\bigg(e^{-x}\sum_{k=0}^{\infty}x^{k}\frac{E[(S_{k+r})_{n,\lambda}]}{k!}	\bigg)\frac{t^{n}}{n!}.\nonumber
\end{align}
Therefore, by \eqref{21}, we obtain the following Dobinski-like formula.
\begin{theorem}
For $n\ge 0$, we have
\begin{displaymath}
\mathrm{Bel}_{n,\lambda}^{(r,Y)}(x)=e^{-x}\sum_{k=0}^{\infty}x^{k}\frac{E[(S_{k+r})_{n,\lambda}]}{k!}.
\end{displaymath}
\end{theorem}
From \eqref{11}, we note that
\begin{align}
	&\frac{1}{m!}E\big[e_{\lambda}^{S_{r}}(t)\big]\Big(E\big[e_{\lambda}^{Y}(t)\big]-1\Big)^{m}\frac{1}{k!}\Big(E\big[e_{\lambda}^{Y}(t)\big]-1\Big)^{k}\label{22}\\
	&=\frac{1}{m!k!}E\big[e_{\lambda}^{S_{r}}(t)\big]\Big(E\big[e_{\lambda}^{Y}(t)\big]-1\Big)^{m+k}\nonumber\\
	&=\frac{(m+k)!}{m!k!}\frac{1}{(m+k)!}\Big(E\big[e_{\lambda}^{Y}(t)\big]-1\Big)^{m+k}E\big[e_{\lambda}^{S_{r}}(t)\big] \nonumber \\
	&=\binom{m+k}{m}\sum_{n=m+k}^{\infty}S_{2,\lambda}^{(r,Y)}(n+r,m+k+r)\frac{t^{n}}{n!}.\nonumber
\end{align}
On the other hand, by \eqref{7} and \eqref{11}, we get
\begin{align}
&\frac{1}{m!}\Big(E\big[e_{\lambda}^{S_{r}}(t)\big]\Big)\Big(E\big[e_{\lambda}^{Y}(t)\big]-1\Big)^{m}\frac{1}{k!}\Big(E\big[e_{\lambda}^{Y}(t)\big]-1\Big)^{k}\label{23} \\
&=\sum_{l=m}^{\infty}S_{2,\lambda}^{(r,Y)}(l+r,m+r)\frac{t^{l}}{l!}\sum_{j=k}^{\infty}S_{2,\lambda}^{Y}(j,k)\frac{t^{j}}{j!}\nonumber \\
&=\sum_{n=k+m}^{\infty}\bigg(\sum_{l=m}^{n-k}\binom{n}{l}S_{2,\lambda}^{(r,Y)}(l+r,m+r)S_{2,\lambda}^{Y}(n-l,k)\bigg)\frac{t^{n}}{n!}.\nonumber
\end{align}
Therefore, by \eqref{22} and \eqref{23}, we obtain the following theorem.
\begin{theorem}
For $m,k\ge 0$ with $n\ge m+k$, we have
\begin{displaymath}
\binom{m+k}{m}S_{2,\lambda}^{(r,Y)}(n+r,m+k+r)=\sum_{l=m}^{n-k}\binom{n}{l}S_{2,\lambda}^{(r,Y)}(l+r,m+r)S_{2,\lambda}^{Y}(n-l,k).
\end{displaymath}
\end{theorem}
From \eqref{15}, \eqref{2-1} and \eqref{7}, we observe that
\begin{align}
&\sum_{n=0}^{\infty}\bigg(\sum_{k=0}^{n}(x)_{k}S_{2,\lambda}^{(r,Y)}(n+r,k+r)\bigg)\frac{t^{n}}{n!}=e^{(x+r)\log E[e_{\lambda}^{Y}(t)]}\label{24} \\
&=\sum_{i=0}^{\infty}(x+r)^{i}\frac{1}{i!}\Big(\log\big(E\big[e_{\lambda}^{Y}(t)\big]-1+1\big)\Big)^{i} \nonumber \\
&=\sum_{i=0}^{\infty}(x+r)^{i}\sum_{j=i}^{\infty}S_{1}(j,i)\frac{1}{j!}\Big(E\big[e_{\lambda}^{Y}(t)\big]-1\Big)^{j}\nonumber \\
&=\sum_{i=0}^{\infty}(x+r)^{i}\sum_{j=i}^{\infty}S_{1}(j,i)\sum_{n=j}^{\infty}S_{2,\lambda}^{Y}(n,j)\frac{t^{n}}{n!}\nonumber \\
&=\sum_{i=0}^{\infty}(x+r)^{i}\sum_{n=i}^{\infty}\sum_{j=i}^{n}S_{1}(j,i)S_{2,\lambda}^{Y}(n,i)\frac{t^{n}}{n!}\nonumber \\
&=\sum_{n=0}^{\infty}\bigg(\sum_{i=0}^{n}\sum_{j=i}^{n}(x+r)^{i}S_{1}(j,i)S_{2,\lambda}^{Y}(n,i)\bigg)\frac{t^{n}}{n!}.\nonumber
\end{align}
Thus, by \eqref{24}, we get the following theorem.
\begin{theorem}
For $n \ge 0$, we have
\begin{displaymath}
\sum_{k=0}^{n}S_{2,\lambda}^{(r,Y)}(n+r,k+r)(x)_{k}= \sum_{i=0}^{n}\sum_{j=i}^{n}S_{1}(j,i)S_{2,\lambda}^{Y}(n,i)(x+r)^{i}.
\end{displaymath}
\end{theorem}

\section{Conclusion}
In recent years, various degenerate versions (see [12,16,23,25,26,30-32]), $\lambda$-analogues (see [15,19,28]) and probabilistic extensions (see [3,10,20-22,29,35]) of many special numbers and polynomials have been explored. Along the way, a lot of interesting results have been found, among which are the degenerate gamma functions (see [24]) and the degenerate umbral calculus (see [18,27]).
In this paper, we considered the probabilistic degenerate $r$-Stirling numbers of the second kind associated with $Y$ and the probabilistic degenerate $r$-Bell polynomials associated with $Y$. They are respectively probabilistic extensions of degenerate versions of the $r$-Stirling numbers of the second kind $S_{2}^{(r)}(n,k)$ and the $r$-Bell polynomials $\mathrm{Bel}_{n}^{(r)}(x)=\sum_{k=0}^{n}S_{2}^{(r)}(n,k)x^{k}$. We obtained by using generating functions their properties, related identities and explicit expressions. \par
As one of our future research projects, we would like to continue to work on probabilistic extensions of some special numbers and polynomials and to find their applications to physics, science and engineering as well as to mathematics.

\end{document}